\documentclass[oneside,reqno]{amsart}

\usepackage{amssymb,latexsym,amsmath,extarrows}
\usepackage{graphicx}

\usepackage{mathtools}

\DeclarePairedDelimiter\floor{\lfloor}{\rfloor}

\newtheorem*{acknowledgement}{Acknowledgements}
\newtheorem{theorem}{Theorem}[section]
\newtheorem{lemma}[theorem]{Lemma}

\newtheorem{remark}[theorem]{Remark}

\newtheorem{corollary}[theorem]{Corollary}

\newtheorem*{definition}{Definition}

\newcommand{\al}{\alpha}
\newcommand{\be}{\beta}

\newcommand{\Ga}{\Gamma}

\newcommand{\De}{\Delta}
\newcommand{\e}{\varepsilon}

\newcommand{\ka}{\kappa}

\newcommand{\La}{\Lambda}
\newcommand{\si}{\sigma}

\newcommand{\om}{\omega}
\newcommand{\Om}{\Omega}

\newcommand{\wh}{\widehat}

\newcommand{\ZR}{\mathbb{R}}
\newcommand{\ZZ}{\mathbb{Z}}

\newcommand{\ZS}{\mathbb{S}}
\newcommand{\ZP}{\mathbb{P}}

\newcommand{\ti}{\tilde}

\newcommand{\hichi}{\raisebox{0.7ex}{\(\chi\)}}

\begin{document}

\title[Fourier decay rates]{Upper bounds for Fourier decay rates of fractal measures}

\author{Xiumin Du}
\address{
Department of Mathematics,
University of Maryland\\
College Park, MD}
\email{xdu@math.umd.edu}

\begin{abstract}
For spherical and parabolic averages of the Fourier transform of fractal measures, we obtain new upper bounds on rates of decay by an ``intermediate dimension'' trick.
\end{abstract}

\maketitle

\section{Introduction}
\setcounter{equation}0
This paper is concerned with average decay rates of the Fourier transform of fractal measures. First recall the notation of ``$\al$-dimensional'' \cite{LR}.

\begin{definition}
Let $0<\alpha\leq d$. We say that $\mu$ is (at least) $\al$-dimensional if it is a positive Borel measure, supported in the unit ball $B^d(0,1)$, that satisfies
$$
c_\al(\mu):=\sup_{x\in\ZR^d, r>0} \frac{\mu(B(x,r))}{r^\al}<\infty.
$$
\end{definition}

Let $S$ be a bounded hypersurface in $\ZR^d$ with everywhere non-vanishing Gaussian curvature and let $d\sigma$ be the induced Lebesgue measure on $S$. We use $\be_d(\al, S)$ to denote the average Fourier decay rate of fractal measures, which is defined as the supremum of the numbers $\be$ for which
\begin{equation} \label{eq:AvrDec}
\left\|\widehat \mu (R\cdot\,)\right\|_{L^2(S, d\sigma)}^2 \lesssim c_\al(\mu)\|\mu\| R^{-\be}
\end{equation}
whenever $R> 1$ and $\mu$ is $\al$-dimensional. In this paper, we will focus on the case $S$ is the unit sphere $\ZS^{d-1}$ or the truncated paraboloid $\ZP^{d-1}$.

The problem of identifying the value of $\beta_d(\al,\ZS^{d-1})$ was proposed by Mattila \cite{M04}, and it relates to the classical distance set conjecture of Falconer \cite{F}.

In dimension two, the exact decay rates are known:
$$
\be_2(\al,S)=
\begin{cases}
\al, & \al\in(0,\,1/2], \quad \text{(Mattila \cite{M})}\\
1/2, & \al\in[1/2,\,1], \quad \text{(Mattila \cite{M})}\\
\al/2, & \al\in[1,\,2], \quad \text{(Wolff \cite{W99})}.
\end{cases}
$$

In higher dimensions, it is known that $\be_d(\al, S) = \alpha$ in the range $\alpha \in (0, \frac{d-1}{2})$, but $\be_d(\al, S)$ is still a mystery for $\frac{d-1}{2}<\al<d$. The current best lower bounds are 
$$
\be_d(\al, S)\geq
\begin{cases}
\al, & \al\in(0,\,\frac{d-1}{2}], \quad \text{(Mattila \cite{M})}\\
\frac{d-1}{2}, & \al\in[\frac{d-1}{2},\,\frac d 2], \quad \text{(Mattila \cite{M})}\\
\frac{(d-1)\al}{d}, & \al\in[\frac d 2,\,d], \quad (\text{D. \emph{et al.} \cite[$d=3$]{DGOWWZ}, D.-Zhang \cite[$d\geq 4$]{DZ}}).
\end{cases}
$$

We remark that the above results were originally computed for either $\ZS^{d-1}$ or $\ZP^{d-1}$. It is however implicit in the arguments given in \cite{M, W99, DGOWWZ, DZ} that the same estimates hold for any bounded hypersurface $S$ with everywhere non-vanishing Gaussian curvature (see, e.g., \cite{CK} for a generalization of \cite{DZ} to a class of hypersurfaces).

Unlike the results for lower bounds, the upper bounds for decay rates are usually obtained by constructing explicit examples and thus the results depend on the hypersurface $S$. The previous best results before this paper are summarized as follows: for the unit sphere, when $d=3$,
$$
\be_3(\al,\ZS^{2}) \leq 
\begin{cases}
\al, & \al\in(0,\,1], \quad \text{(\cite[Chapter 15.2]{M'})}\\
\frac{\al+1}{2}, & \al\in[1,\,3], \quad (\text{Knapp example}),
\end{cases}
$$
and when $d\geq 4$,
$$
\be_d(\al,\ZS^{d-1}) \leq 
\begin{cases}
\al, & \al\in(0,\,\frac d 2], \quad \text{(\cite{M'})}\\
\al-1+\frac{2(d-\al)}{d}, & \al\in[\frac d 2,\,d], \quad (\text{Luc\`a-Rogers \cite{LR}});
\end{cases}
$$
for the truncated paraboloid and $d\geq 3$,
$$
\be_d(\al,\ZP^{d-1}) \leq 
\begin{cases}
\al, & \al\in(0,\,\frac{d-1}{2}], \quad \text{(\cite{M'})}\\
\frac{(d-1)(\al+1)}{d+1}, & \al\in[\frac{d-1}{2},\,d], \quad (\text{Barcel\'o \emph{et al.} \cite{BBCRV}}).
\end{cases}
$$

It is worth mentioning that when $\al=d-1$, one can find a better upper bound of $\frac{(d-1)^2}{d}$ by examining an example of Bourgain \cite{jB16} carefully. As this upper bound coincides with the lower bound established in \cite{DGOWWZ, DZ}, the exact decay rate can be determined in this case: $$\be_d(d-1,\ZP^{d-1}) = \frac{(d-1)^2}{d}.$$ Bourgain's example is a Schr\"odinger solution essentially supported in a small neighborhood of a hyperplane. 
%Being connected with many topics, it is natural and interesting to explore the full strength of Bourgain's idea. to all intermediate dimensions. 
Recently, the authors of \cite{DKWZ} extended Bourgain's idea to intermediate dimensions and disproved Schr\"odinger maximal estimates in certain range. In this paper, we further explore this ``intermediate dimension'' trick to adapt the examples from \cite{BBCRV, LR} and obtain improved upper bounds of Fourier decay rates.

We first state the results for spheres. For convenience of notation, we introduce the following functions $\ka_1$ and $\ka_2$:
\begin{equation}\label{k1}
    \ka_1(m;\al,d):=\frac{d-m/2-\al}{d-m}, \quad \ka_2(m;\al,d):=\frac{d-\al}{2(d-m)}\,.
\end{equation}
For $\ka_1$ and $\ka_2$, we are only interested in the cases that $\al\in(d/2,d)$ and $m$ is an integer with $0<m<d/2$. In this range, for fixed $\al$ and $d$, as $m$ increases, $\ka_1(m;\al,d)$ decreases and $\ka_2(m;d,\al)$ increases.

\begin{theorem} \label{thm:sph}
Let $d\geq 4$ and $\al\in (d/2,d)$. Then
$$
\be_d(\al,\ZS^{d-1})\leq \al-1+2\ka(\al,d)\,,
$$
where $\ka(\al,d)$ is given as follows:

\noindent (a). For $\al\in[d-1,d)$,
$$
\ka(\al,d)=\ka_2(1;\al,d)=\frac{d-\al}{2(d-1)}\,;
$$

\noindent (b). For $\al \in [d-j,d-j+1]$ with $j=2,3,\cdots,\floor{\frac{d-1}{2}}$,
\begin{align*}
    \ka(\al,d)&=\min\Big\{\ka_1(j-1;\al,d),\,\,\ka_2(j;\al,d)\Big\}\\
    &=
    \begin{cases}
    \ka_2(j;\al,d), & d-j \leq \al \leq d-j+\frac{d-2j}{d-j-1}\,,\vspace{0.1 in}\\
    \ka_1(j-1;\al,d), & d-j+\frac{d-2j}{d-j-1}\leq \al \leq d-j+1 \,;
    \end{cases}
\end{align*}

\noindent (c). For $d$ even and $\al\in(\frac d 2,\frac d2+1]$,
$$
\ka(\al,d)=\ka_1\left(\frac d2-1;\al,d\right)=\frac{3d+2-4\al}{2(d+2)}\,;
$$

\noindent (d). For $d$ odd and $\al\in(\frac d 2,\frac {d+1}{2}]$,
$$
\ka(\al,d)=\ka_1\left(\frac{d-1}{2};\al,d\right)=\frac{3d+1-4\al}{2(d+1)}\,.
$$
\end{theorem}
\bigskip
%\begin{remark}
Note that the previous best result from  \cite{LR} is equivalent to saying that for $d\geq 4$ and $ \al\in(d/2,d)$,
$$
\be_d(\al,\ZS^{d-1})\leq \al-1+2\ka_1(0;\al,d)\,.
$$
Since $\ka_1(m;\al,d)$ is a decreasing function of $m$ and 
$$
\ka_2(m;\al,d)< \ka_1(0;\al,d) \quad \text{for} \quad m<\frac d 2, 
$$
we see that Theorem \ref{thm:sph} is indeed better in the whole range stated in the theorem.
%\end{remark}

Next, we turn to the paraboloids. Define three more functions:
\begin{equation}\label{k3}
\begin{split}
    \ka_3(m;\al,d):=\frac{d-m/2-\al}{d-m+1}, \quad &\ka_4(m;\al,d):=\frac{d-\al}{2(d-m+1)},\\
    \ka_5(m;\al,d):=&\frac{d-\al-1}{2(d-m-1)}\,.
\end{split}
\end{equation}
Here $m$ is again a positive integer. For $\ka_3$, we will focus on the range $\al\in(\frac{d-1}{2},d)$ and  $0<m<d/2$; for $\ka_4$, consider the cases $\al\in(\frac{d-1}{2},d)$ and  $0<m\leq d/2$; for $\ka_5$, we are interested in the situation that $\al\in(\frac{d-1}{2},d-1)$ and $0<m\leq d/2$. In all these cases, for fixed $\al$ and $d$, as $m$ increases, $\ka_3(m;\al,d)$ decreases, $\ka_4(m;d,\al)$ and $\ka_5(m;d,\al)$ increase.

\begin{theorem} \label{thm:prb}
Let $d\geq 3$ and $\al\in (\frac{d-1}{2},d)$. Then
$$
\be_d(\al,\ZP^{d-1})\leq \al-1+2\ti\ka(\al,d)\,,
$$
where $\ti\ka(\al,d)$ is given as follows:

\noindent (a). For $\al\in[d-1,d)$,
$$
\ti\ka(\al,d)=\ka_4(1;\al,d)=\frac{d-\al}{2d}\,;
$$

\noindent (b). For $\al \in [d-j,d-j+1]$ with $2\leq j \leq \floor{\frac{d+1}{3}}$,
\begin{align*}
    \ti\ka(\al,d)&=\min\Big\{\ka_3(j-1;\al,d),\,\,\ka_4(j;\al,d)\Big\}\\
    &=
    \begin{cases}
    \ka_4(j;\al,d), & d-j \leq \al \leq d-j+\frac{d-2j+1}{d-j}\,,
    \vspace{0.1 in}
    \\
    \ka_3(j-1;\al,d), & d-j+\frac{d-2j+1}{d-j}\leq \al \leq d-j+1 \,;
    \end{cases}
\end{align*}

\noindent (c). For $\al \in [d-j,d-j+1]$ with $j = \floor{\frac{d+1}{3}} +1$,
$$
\ti\ka(\al,d)=\ka_3\left(\floor{\frac{d+1}{3}};\al,d\right) \,;
$$

\noindent (d). For $\al \in [d-j,d-j+1]$ with $\floor{\frac{d+1}{3}}+2 \leq j \leq \floor{\frac d2}$,
\begin{align*}
    &\ti\ka(\al,d)\\
    =&\min\Big\{\ka_3(j-2;\al,d),\,\,\max\big\{\ka_3(j-1;\al,d),\,\,\ka_5(j-1;\al,d)\big\},\,\,\ka_5(j;\al,d)\Big\}\\
    =&
    \begin{cases}
    \min\Big\{\ka_3(j-1;\al,d),\,\,\ka_5(j;\al,d)\Big\}, & d-j \leq \al \leq d-j+\frac{2(d-2j+1)}{d-j-2}\,,
    \vspace{0.1 in} \\
    \min\Big\{\ka_3(j-2;\al,d),\,\,\ka_5(j-1;\al,d)\Big\}, & d-j+\frac{2(d-2j+1)}{d-j-2}\leq \al \leq d-j+1 \,;
    \end{cases}
\end{align*}

\noindent (e). For $d$ odd, $d\geq 7$ and $\al\in(\frac{d-1}{2},\frac {d+1}{2}]$,
$$
\ti\ka(\al,d)=\ka_3\left(\frac{d-3}{2}\right)=\frac{3d+3-4\al}{2(d+5)}\,;
$$
Note that the cases $d=3,5$ and $\al\in(\frac{d-1}{2},\frac {d+1}{2}]$ were covered in part (c).

\noindent (f). For $d$ even and $\al\in(\frac {d-1}{2},\frac {d}{2}]$,
$$
\ti\ka(\al,d)=\ka_3\left(\frac{d}{2}-1;\al,d\right)=\frac{3d+2-4\al}{2(d+4)}\,.
$$
\end{theorem}
\bigskip
%\begin{remark}
Note that the previous best upper bound from  \cite{BBCRV} is equivalent to saying that for $d\geq 3$ and $ \al\in(\frac{d-1}{2},d)$,
$$
\be_d(\al,\ZP^{d-1})\leq \frac{(d-1)(\al+1)}{d+1}= \al-1+2\ka_3(0;\al,d)\,.
$$
Since $\ka_3(m;\al,d)$ is a decreasing function of $m$ and 
$$
\ka_4(m;\al,d)< \ka_3(0;\al,d) \quad \text{for} \quad m<\frac{d+1}{2}, 
$$
we see that Theorem \ref{thm:prb} is an improvement in the whole range stated in the theorem.
%\end{remark}

\begin{remark}
It is straightforward to check $\ti\ka(\al,d)<\ka(\al,d)$. In other words, the examples for parabolic decay rates are better than those for spherical decay rates.
\end{remark}

%\begin{remark}
By combining part (a) of Theorem \ref{thm:prb} and the lower bounds from \cite{DGOWWZ, DZ}, we can now determine the exact value of the parabolic Fourier decay rates for $\al\in [d-1,d)$. We record this result in the following corollary.
\begin{corollary}
Let $d-1\leq \al<d$ and $d\geq3$. Then
\begin{equation}
    \be_d(\al,\ZP^{d-1})= \al-1+\frac{d-\al}{d}= \frac{(d-1)\al}{d}
\end{equation}
\end{corollary}

%\end{remark}

\begin{remark}
To get a feeling about the numerology in Theorem \ref{thm:prb}, let's explicitly write out $\ti\ka(\al,d)$ with $\al\in(\frac{d-1}{2},d-1]$ for some small values of $d$. This will also be useful in the next remark.
\begin{itemize}
    \item For $d=3,4$, $$\ti\ka(\al,d)=\ka_3(1;\al,d)=\frac{2d-1-2\al}{2d}, \quad \frac{d-1}{2}<\al\leq d-1\,;$$
    \item For $d=5,6,7$,
    $$\ti\ka(\al,d)=
    \begin{cases}
    \ka_3(2;\al,d)=\frac{d-1-\al}{d-1}, \quad &\frac{d-1}{2}<\al\leq d-2, \vspace{0.1 in}\\
    \ka_4(2;\al,d)=\frac{d-\al}{2(d-1)},&d-2\leq \al\leq d-2+\frac{d-3}{d-2}, \vspace{0.1 in}\\
    \ka_3(1;\al,d)=\frac{2d-1-2\al}{2d},&d-2+\frac{d-3}{d-2}\leq \al\leq d-1\,;
    \end{cases}
    $$
\end{itemize}
The situation becomes more complicated for larger $d$, and $\ka_5(m;\al,d)$ will also come into play when $d$ is large enough.
\end{remark}

\begin{remark} Let us see what we can tell about Falconer's distance set conjecture from our new theorems.

(a). For $\al$ close to and greater than $d/2$, Theorem \ref{thm:prb} tells us that $\be_d(\al,\ZP^{d-1})\leq \al-1+2\ti\ka(\al,d)$, where
$$
\ti\ka(\al,3)=\ka_3(1;\al,3)=\frac{5-2\al}{6}, \quad \ti\ka(\al,5)=\ka_3(2;\al,5)=\frac{4-\al}{4},
$$
$$
\ti\ka(\al,d)=\ka_3\left(\frac{d-3}{2};\al,d\right)=\frac{3d+3-4\al}{2(d+5)}\quad \text{ for } d \text{ odd and } d\geq 7,
$$
and
$$
\ti\ka(\al,d)=\ka_3\left(\frac{d}{2}-1;\al,d\right)=\frac{3d+2-4\al}{2(d+4)}\quad \text{ for } d \text{ even and } d\geq 4.
$$

(b). According to a famous scheme developed by Mattila, the Fourier decay rates of fractal measures and Falconer's conjecture are related as follows (see for example \cite{DGOWWZ}): 

Suppose that \eqref{eq:AvrDec} holds for $S=\ZS^{d-1}$ with some $\be\geq d-\al$. Then Falconer's distance set conjecture holds for $\al$, i.e.  for any compact subset $E$ of $\ZR^d$,
$$
\rm{dim} (E)>\al \implies |\De(E)|>0,
$$
where $|\cdot|$  denotes the Lebesgue measure, $\text{dim}(\cdot)$ is the Hausdorff dimension and $\De(E)$ is the distance set given by $\De(E)=\left\{|x-y|:x,y\in E\right\}\,.$  
The threshold for $\al$ in Falconer's conjecture is $d/2$. 

(c). Suppose we plan to approach Falconer's conjecture using the above relation. Assume \eqref{eq:AvrDec} also holds for $S=\ZP^{d-1}$ with the same $\be\geq d-\al$. (This is the case in all previous works \cite{M, W99, DGOWWZ, DZ}). Then Theorem \ref{thm:prb} tells us that the best possible threshold for $\al$ one could get using Mattila's scheme is
$$
\frac 7 4=\frac 32+\frac 14 \text { when } d=3,\quad \frac 8 3=\frac 52+\frac 16 \text { when } d=5,
$$
$$
\frac d2+\frac{1}{d+3} \text { when $d$ odd and } d\geq 7, \quad \frac d2+\frac{1}{d+2} \text { when $d$ even and } d\geq 4.
$$

This suggests that new approach (e.g., \cite{bL,GIOW}) may be needed to fully resolve Falconer's conjecture.
\end{remark}

\vspace{.1in}

\noindent \textbf{Notation.} We write $A\lesssim B$ if $A\leq CB$ for some absolute constant $C$, $A \sim B$ if $A\lesssim B$ and $B\lesssim A$, and $A\lessapprox B$ if $A\leq C_\e R^\e B$ for any $\e>0, R>1$. Let $c=1/1000$ be fixed. By $\rho$-lattice points in $\ZR^d$ we mean the points in $\rho \ZZ^d$. Let $B^d(x,r)$ denote the ball centered at $x$, of radius $r$, in $\ZR^d$. 

\begin{acknowledgement}
 The author is supported by the National Science Foundation under Grant No. DMS-1856475.
\end{acknowledgement}

\section{Proof of Theorem \ref{thm:sph} - Spherical decay rates}\label{sec-sph}

Let $\mu$ be $\alpha$-dimensional. Given a function $g$ on the unit ball $B^d(0,1)$, we can write $g=g_1-g_2+i(g_3-g_4)$, where each component $g_j$ is positive. Then by considering the positive measures $g_j\mu$, the estimate \eqref{eq:AvrDec} tells us that
$$
\left\|\wh{g\mu}(R\cdot\,)\right\|^2_{L^2(S)}\lesssim c_\al(\mu)\|\mu\|R^{-\be} \|g\|^2_{L^\infty}\,.
$$
Thus, by duality, we are looking for an upper bound for the $\be$ such that
\begin{equation} \label{eq:L1} 
\left\|E_Sf(R\cdot\,)\right\|_{L^1(d\mu)}\lesssim R^{-\be/2}\sqrt{c_\al(\mu)\|\mu\|} \|f\|_{L^2(S)}\,,
\end{equation}
where 
$$
E_Sf(x)=(fd\si)^\vee(x)=\frac{1}{(2\pi)^{d/2}}\int_S e^{i\omega\cdot x}f(\omega)\,d\si(\om)\,.
$$

\vspace{.25 in}

This example is adapted from that of \cite{LR}.
Let $c=1/1000$ be a fixed small constant and $0<\ka<1/2$. The exact value of $\ka$ will be chosen later.
Let  $1\leq m< d/2$ and $d\geq 4$. Denote
$$x=(x_1,\cdots,x_d)=(x',x'')\in B^d(0,1)\,,$$ $$\xi=(\xi_1,\cdots,\xi_d)=(\xi',\xi'')\in \ZS^{d-1}\,,$$
where $$
x'=(x_1,\cdots,x_m), \quad
x''=(x_{m+1},\cdots,x_d),$$
$$
\xi'=(\xi_1,\cdots,\xi_m), \quad
\xi''=(\xi_{m+1},\cdots,\xi_d).$$

\noindent For $S=\ZS^{d-1}$, the unit sphere in $\ZR^d$, we write $E_Sf(Rx)$ as
\begin{equation}\label{Ef'}
   Ef(Rx)=\frac{1}{(2\pi)^{d/2}} \int_{\ZS^{d-1}} e^{i(Rx'\cdot \xi'+Rx''\cdot \xi'')}f(\xi)\,d\si(\xi). 
\end{equation}

To prove Theorem \ref{thm:sph}, we'll test the estimate \eqref{eq:L1} on the characteristic function $f(\xi)=\hichi_\Om(\xi)$, where the set $\Om$ is defined by
\begin{equation} \label{Om'}
    \Om:=\left[B^m(0,cR^{-1/2})\times \left(\Ga+B^{d-m}(0,cR^{-1})\right)\right]\cap \ZS^{d-1}\,,
\end{equation}
and
\begin{equation}\label{Ga}
    \Ga:=\left\{\om\in \ZS^{d-m-1}\,:\, R^\ka \om\in 2\pi\ZZ^{d-m}\right\}\,.
\end{equation}
So we have that $\|f\|_2=\si(\Om)^{1/2}$.

It's well known (see, for example, a survey about lattice points on spheres \cite{fF}) that for $d-m\geq 2$, there holds
$$
\#\Ga \gtrapprox R^{\ka(d-m-2)}\,,
$$
for a sequence of $R$ tending to $\infty$. We'll focus on such values of $R$.
Note that, in the definition of $\Om$, each point in $\Ga$ gives us a small patch on $\ZS^{d-1}$, which has size $\sim R^{-1/2}$ in $m$ dimension and $\sim R^{-1}$ in each of the other $(d-m-1)$ dimensions. Therefore,
\begin{equation} \label{sizeOm}
 \si(\Om) \gtrapprox R^{\ka(d-m-2)-\frac m2-(d-m-1)}=R^{\ka(d-m-2)-d+\frac m2+1}\,.
\end{equation}

Next, we define a set $\La$ in $B^d(0,1)$ by \begin{equation}\label{La''}
    \La:=\left[B^m(0,cR^{-1/2})\times \left(R^{\ka-1}\ZZ^{d-m}+B^{d-m}(0,cR^{-1})\right)\right] \cap B^d(0,1)\,.
\end{equation}
The idea is that for $x\in \La$, the phase of the integrand in \eqref{Ef'} is sufficiently close to $2\pi i \ZZ$, and so there is little cancellation - see Lemma \ref{C:sizeEf'}. Now define $\mu$ by
\begin{equation} \label{mu'}
    d\mu=\hichi_\La dx,
\end{equation}
where $dx$ is the Lebesgue measure in $\ZR^d$. From the definition it follows that
\begin{equation}\label{sizeLa'}
    \|\mu\|=|\La| \quad \text{and} \quad |\La|\sim R^{-m/2}\left(R^{1-\ka}R^{-1}\right)^{d-m}=R^{-\ka(d-m)-m/2}\,.
\end{equation}

We need the following two lemmas, whose proofs are postponed.

\begin{lemma}\label{C:sizeEf'}
For $f$ given above, 
\begin{equation} \label{sizeEf'}
 |Ef(Rx)|\sim \si(\Om), \quad \quad \forall x\in \La\,.  
\end{equation}
\end{lemma}

\begin{lemma} \label{C:si'}
By taking 
\begin{equation}\label{ka}
    \ka=\begin{cases} 
    \ka_1(m;\al,d)=\frac{d-m/2-\al}{d-m},& \al \in (\frac d 2, d-m] \vspace{0.1 in}\\
    \ka_2(m;\al,d)=\frac{d-\al}{2(d-m)}, & \al\in [d-m,d),
    \end{cases}
\end{equation}
we have
\begin{equation}\label{c_al'}
    c_\al(\mu) \sim R^{\al-d}\,.
\end{equation}
\end{lemma}

By plugging in \eqref{sizeOm}, \eqref{sizeLa'}, \eqref{sizeEf'} and \eqref{c_al'}, we obtain
$$
\frac{\|Ef(R\cdot\,)\|_{L^1(d\mu)}}{\sqrt{c_\al(\mu)\|\mu\|}\|f\|_2} \sim \frac{\si(\Om)|\La|}{R^{(\al-d)/2}|\La|^{1/2}\si(\Om)^{1/2}}\gtrapprox R^{-\ka+\frac{1-\al}{2}}\,.
$$
Comparing the above with \eqref{eq:L1}, letting $R$ tend to infinity and taking $\be$ sufficiently close to $\be_{d}(\al,\ZS^{d-1})$, we see that 
$$
\be_{d}(\al,\ZS^{d-1}) \leq\al-1+2\ka\,,
$$
where $\ka$ is given as in \eqref{ka}. To prove Theorem \ref{thm:sph}, we just take suitable $m$ for different values of $\al$. It follows directly from \eqref{ka} that we can choose $\ka$ as follows:
\begin{itemize}
    \item For $\al \in [d-1,d)$, $\ka=\ka_2(1;\al,d)$.
    \vspace{0.1 in}
    \item For $d$ even and $\al\in(\frac d 2,\frac d2+1]$, $\ka=\ka_1\left(\frac d2-1;\al,d\right)$.
    \vspace{0.1 in}
    \item For $d$ odd and $\al\in(\frac d 2,\frac {d+1}{2}]$,
$\ka=\ka_1\left(\frac{d-1}{2};\al,d\right).$
\vspace{0.1 in}
\item For $\al \in [d-j,d-j+1]$ with $j=2,3,\cdots,\floor{\frac{d-1}{2}}$,
$$
\ka=\min\Big\{\ka_1(j-1;\al,d),\,\,\ka_2(j;\al,d)\Big\}\,.
$$
It is straightforward to check that 
$$
\ka_2(j;\al,d) \leq \ka_1(j-1;\al,d) \quad \iff \quad \al \leq  d-j+\frac{d-2j}{d-j-1}\,.
$$
Also note that $0<\frac{d-2j}{d-j-1}<1$ in this case.
\end{itemize}

This finishes the proof of Theorem \ref{thm:sph} up to Lemmas \ref{C:sizeEf'} and \ref{C:si'}.

\subsection{Proof of Lemma \ref{C:sizeEf'}}
Since $f=\hichi_\Om$, we have
$$
Ef(Rx)=\frac{1}{(2\pi)^{d/2}} \int_{\Om} e^{i(Rx'\cdot\xi'+Rx''\cdot \xi'')} \,d\si(\xi). 
$$
So it suffices to prove that 
\begin{equation} \label{phase'}
 Rx'\cdot\xi'+R x''\cdot \xi''\in 2\pi \ZZ +(-\frac{1}{100},\frac{1}{100})\,,
\end{equation}
provided that $\xi\in \Om$ and $x\in\La$. Indeed, by definitions of $\Om$ and $\La$, we write
$$
|\xi'|<cR^{-\frac 12}, \quad |x'|<cR^{-\frac 1 2}
$$
$$
\xi''=2\pi R^{-\ka}m+v,\quad \text{where} \quad m\in\ZZ^{d-m}, |m|<\frac{1}{2\pi}R^{\ka}, |v|<cR^{-1},
$$
and
$$
x''=R^{\ka-1}\ell+u,\quad \text{where} \quad \ell\in\ZZ^{d-m}, |\ell|< R^{1-\ka}, |u|<cR^{-1}\,.
$$
Then it is straightforward to verify that \eqref{phase'} holds.
\begin{itemize}
    \item $|Rx'\cdot \xi'|<RcR^{-1/2}cR^{-1/2}=c^2$.
    \vspace{0.1 in}
    \item For $Rx''\cdot \xi''$, we have
   \begin{align*}
       Rx''\cdot \xi''=&R(R^{\ka-1}\ell+u)\cdot(2\pi R^{-\ka}m+v)\\
        =&2\pi \ell\cdot m+R^{\ka}\ell\cdot v +2\pi R^{1-\ka}u\cdot m +Ru\cdot v\,,
   \end{align*}
   where $2\pi\ell\cdot m\in2\pi\ZZ$ and the other three terms are bounded by
   $$ R^\ka R^{1-\ka} cR^{-1}+R^{1-\ka}cR^{-1}R^\ka+RcR^{-1}cR^{-1}=
   c+c+c^2R^{-1}\,.$$
\end{itemize}
Therefore, \eqref{phase'} follows by taking $c$ sufficiently small, say $c=1/1000$.

\subsection{Proof of Lemma \ref{C:si'}}
Recall that $d\mu=\hichi_\La\,dx$ and $\La$ is defined by
\begin{equation}\label{La'''}
    \La:=\left[B^m(0,cR^{-1/2})\times \left(R^{\ka-1}\ZZ^{d-m}+B^{d-m}(0,cR^{-1})\right)\right] \cap B^d(0,1)\,.
\end{equation}
We aim to prove that
$$
c_\al(\mu) \sim R^{\al-d}\,,
$$
by taking 
\begin{equation}\label{ka'}
    \ka=\begin{cases} 
    \ka_1(m;\al,d)=\frac{d-m/2-\al}{d-m},& \al \in (\frac d 2, d-m] \vspace{0.1 in} \\
    \ka_2(m;\al,d)=\frac{d-\al}{2(d-m)}, & \al\in [d-m,d).
    \end{cases}
\end{equation}

For convenience, we write
$$
c_\al(\mu):=\sup_{x\in\ZR^d, r>0} \frac{\mu(B(x,r))}{r^\al}=\sup_{r>0}c_\al(\mu,r)\,,
$$
where
$$
c_\al(\mu,r):=\sup_{x\in\ZR^d} \frac{\mu(B(x,r))}{r^\al}\,.
$$
We will calculate $c_\al(\mu,r)$ directly from \eqref{La'''}.
The important scales for $r$ are ordered as follows:
$$
R^{-1}<R^{\ka-1}<R^{-1/2}<1\,.
$$
Now we calculate $C_\al(\mu,r)$ for different values of $r$.

\begin{itemize}
    \item For $0<r\leq R^{-1}$,
    $$
    c_\al(\mu,r)\sim \frac{r^d}{r^\al}=r^{d-\al}\,.
    $$
    
    \noindent Since $d-\al>0$, we have
    \begin{equation}\label{1''}
        \sup_{0<r\leq R^{-1}} c_\al(\mu,r)\sim  c_\al(\mu,R^{-1})\sim R^{\al-d}\,.
    \end{equation}
    
    \item For $R^{-1}\leq r\leq R^{\ka-1}$,
    $$
    c_\al(\mu,r)\sim \frac{r^m\cdot R^{-(d-m)}}{r^\al}=r^{m-\al}R^{-(d-m)}\,.
    $$
    
    \noindent Since $m<d/2< \al$, we have
    \begin{equation}\label{2''}
        \sup_{R^{-1}\leq r\leq R^{\ka-1}} c_\al(\mu,r)\sim  c_\al(\mu,R^{-1})\,.
    \end{equation}
    
    \item For $R^{\ka-1}\leq r\leq R^{-\frac 1 2}$,
    $$
    c_\al(\mu,r)\sim \frac{r^m\cdot \left(\frac{r}{R^{\ka-1}}R^{-1}\right)^{d-m}}{r^\al}=r^{d-\al}R^{-\ka (d-m)}\,.
    $$
    
    \noindent Since $d-\al>0$, we have
    \begin{equation}\label{3''}
        \sup_{R^{\ka-1}\leq r\leq R^{- 1/2}} c_\al(\mu,r)\sim  c_\al(\mu,R^{-\frac 1 2})\sim R^{-\frac{d-\al}{2}-\ka(d-m)}\,.
    \end{equation}
    
     \item For $R^{-\frac 1 2}\leq r\leq 1$,
    $$
    c_\al(\mu,r)\sim \frac{R^{-m/2}\cdot \left(\frac{r}{R^{\ka-1}}R^{-1}\right)^{d-m}}{r^\al}=r^{d-m-\al}R^{-\ka (d-m)-m/2}\,.
    $$
     
     \noindent If $\al\leq d-m$, we have
    \begin{equation}\label{4''}
        \sup_{R^{-1/2}\leq r\leq 1} c_\al(\mu,r)\sim  c_\al(\mu,1)\sim R^{-\ka(d-m) -\frac m2}\,,
    \end{equation}
     
     \noindent and if $\al\geq d-m$, we have
    \begin{equation}\label{5''}
        \sup_{R^{-1/2}\leq r\leq 1} c_\al(\mu,r)\sim  c_\al(\mu,R^{-\frac 12})\,.
    \end{equation}
\end{itemize}

\noindent It is also obvious that 
$$
\sup_{r\geq 1} c_\al(\mu,r)\sim  c_\al(\mu,1)\,.
$$
Therefore, for $\al\leq d-m$, by combining \eqref{1''}, \eqref{2''}, \eqref{3''} and \eqref{4''},  we can tell that
\begin{align*}
    c_\al(\mu)&\sim \max\left\{c_\al(\mu, R^{-1}),c_\al(\mu,1)\right\} \\
    &\sim \max\left\{R^{\al-d}, R^{-\ka(d-m)-\frac m2}\right\} =R^{\al-d}\,,
\end{align*}
provided that
$$
\ka=\ka_1(m;\al,d)=\frac{d-m/2-\al}{d-m}\,.
$$

\noindent And for $\al\geq d-m$, by combining \eqref{1''}, \eqref{2''}, \eqref{3''} and \eqref{5''}, we can tell that
\begin{align*}
    c_\al(\mu)&\sim \max\left\{c_\al(\mu, R^{-1}),c_\al(\mu,R^{-\frac 12})\right\} \\
    &\sim \max\left\{R^{\al-d}, R^{-\frac{d-\al}{2}-\ka (d-m)}\right\} =R^{\al-d}\,,
\end{align*}
provided that
$$
 \ka=\ka_2(m;\al,d)=\frac{d-\al}{2(d-m)}\,,
$$
as desired. This completes the proof of Lemma \ref{C:si'}.

\section{Proof of Theorem \ref{thm:prb} - Parabolic decay rates}\label{sec-prb}

This example is adapted from that of \cite{BBCRV} in a similar way as in the previous section.
Recall that $c=1/1000$ is a fixed small constant. In this section, we will still use but redefine the notations $f, \Om$ and $\La$.
Let $0<\ka<1/2$. Let $d\geq 3$ and $1\leq m\leq d/2$. In $\ka_3(m;\al,d)$ below, $m<d/2$, while in  $\ka_4(m;\al,d)$ and $\ka_5(m;\al,d)$ below, $m$ could be $d/2$.
Denote
$$x=(x_1,\cdots,x_d)=(x',x'',x_d)\in B^d(0,1)\,,$$ $$\xi=(\xi_1,\cdots,\xi_{d-1})=(\xi',\xi'')\in B^{d-1}(0,1)\,,$$
where $$
x'=(x_1,\cdots,x_m), \quad
x''=(x_{m+1},\cdots,x_{d-1}),$$
$$
\xi'=(\xi_1,\cdots,\xi_m), \quad
\xi''=(\xi_{m+1},\cdots,\xi_{d-1}).$$

\noindent For $S=\ZP^{d-1}$, the truncated paraboloid in $\ZR^d$, we write $E_Sf(Rx)$ as
\begin{equation}\label{Ef}
   Ef(Rx)=\frac{1}{(2\pi)^{d/2}} \int_{B^{d-1}(0,1)} e^{iR(x'\cdot\xi'+x''\cdot \xi'' +x_d|\xi'|^2+x_d|\xi''|^2)}f(\xi)\,d\xi. 
\end{equation}

For simplicity, we denote $B^d(0,r)$ by $B^d_r$,  and write the interval $(-r,r)$ as $I_r$.
To prove Theorem \ref{thm:prb}, we'll test the estimate \eqref{eq:L1} on the characteristic function $f(\xi)=\hichi_\Om(\xi)$, where the set $\Om$ is defined by
\begin{equation} \label{Om}
    \Om:=\left[B^m_{cR^{-1/2}} \times \left(2\pi R^{-\ka} \ZZ^{d-m-1}+B^{d-m-1}_{cR^{-1}}\right)\right]\cap B^{d-1}(0,1)\,.
\end{equation}
By definition, we have
\begin{equation}\label{OmSize}
   \|f\|_2=|\Om|^{1/2}\quad \text{and} \quad |\Om|\sim R^{(\ka-1)(d-m-1)-m/2}\,.
\end{equation}

Next, we define a set $\La$ in $B^d(0,1)$ by \begin{equation}\label{La}
    \La:=\left[B^m_{cR^{-1/2}}\times \left(R^{\ka-1}\ZZ^{d-m-1}+B^{d-m-1}_{cR^{-1}}\right)\times\left(\frac{1}{2\pi}R^{2\ka-1}\ZZ+I_{cR^{-1}}\right)\right] \cap B^d(0,1)\,,
\end{equation}
Now, define $\mu$ by
\begin{equation} \label{mu}
    d\mu=\hichi_\La dx,
\end{equation}
where $dx$ is the Lebesgue measure in $\ZR^d$. From the definition it follows that
\begin{equation}\label{sizeLa}
    \|\mu\|=|\La| \quad \text{and} \quad |\La|\sim R^{-m/2-\ka(d-m-1)-2\ka}=R^{-\ka(d-m+1)-m/2}\,.
\end{equation}

Moreover, we have the following two lemmas, whose proofs are postponed.

\begin{lemma}\label{C:sizeEf}
For $f$ given above, 
\begin{equation} \label{sizeEf}
 |Ef(Rx)|\sim |\Om|, \quad \quad \forall x\in \La\,.  
\end{equation}
\end{lemma}

\begin{lemma} \label{C:si}
We have
\begin{equation}\label{c_al}
    c_\al(\mu) \sim R^{\al-d}\,,
\end{equation}
by taking $\ka$ as follows:

(a). If $1\leq m\leq \frac{d+1}{3}$, then
\begin{equation} \label{a1}
    \ka=\ka_3(m;\al,d)\quad \text{for } \quad m\leq \al\leq d-m, 
\end{equation}
and
\begin{equation} \label{a2}
    \ka=\ka_4(m;\al,d)\quad \text{for } \quad d-m\leq \al< d. 
\end{equation}

 (b). If $\frac{d+1}{3}<m< \frac d2$, then
\begin{equation} \label{b1}
    \ka=\ka_3(m;\al,d)\quad \text{for } \quad  m\leq \al\leq d-m-1,
\end{equation}
and
\begin{align} \label{b2}
    \ka&=\max\Big\{\ka_3(m;\al,d),\,\,\ka_5(m;\al,d)\Big\} \quad \text{for } \quad d-m-1\leq \al\leq d-m \\
    &=
    \begin{cases}
    \ka_3(m;\al,d), & d-m-1 \leq \al \leq d-m-1+\frac{2(d-2m-1)}{d-m-3}\,,\vspace{0.1 in} \notag\\
    \ka_5(m;\al,d), & d-m-1+\frac{2(d-2m-1)}{d-m-3}\leq \al \leq d-m \,,
    \end{cases}
\end{align}
and
\begin{equation}\label{b3}
    \ka=\ka_5(m;\al,d)\quad \text{for } \quad  d-m\leq \al\leq \frac{d+m-1}{2}, 
\end{equation}
and
\begin{equation}\label{b4}
    \ka=\ka_4(m;\al,d)\quad \text{for } \quad  \frac{d+m-1}{2} \leq \al <d.
\end{equation}
Moreover, \eqref{b3} and \eqref{b4} also holds when $m=\frac d 2$.
\end{lemma}

By plugging in \eqref{OmSize}, \eqref{sizeLa}, \eqref{sizeEf} and \eqref{c_al}, we obtain
$$
\frac{\|Ef(R\cdot\,)\|_{L^1(d\mu)}}{\sqrt{c_\al(\mu)\|\mu\|}\|f\|_2} \sim \frac{|\Om||\La|}{R^{(\al-d)/2}|\La|^{1/2}|\Om|^{1/2}}\sim R^{-\ka+\frac{1-\al}{2}}\,.
$$
Comparing the above with \eqref{eq:L1}, letting $R$ tend to infinity and taking $\be$ sufficiently close to $\be_{d}(\al,\ZP^{d-1})$, we see that 
\begin{equation} \label{beka}
    \be_{d}(\al,\ZP^{d-1}) \leq\al-1+2\ka\,,
\end{equation}
where $\ka$ is given as in Lemma \ref{C:si}. To prove Theorem \ref{thm:prb}, we just take suitable $m$ for different values of $\al$:

\begin{itemize}
    \item For $\al \in [d-1,d)$, by \eqref{a2} we can take 
    $$\ka=\ka_4(1;\al,d)$$.
    
    \item For $\al \in [d-j,d-j+1]$ with $2\leq j \leq \floor{\frac{d+1}{3}}$, by \eqref{a1} we can take $\ka=\ka_3(j-1;\al,d)$, and by \eqref{a2} we can take $\ka=\ka_4(j;\al,d)$. Therefore, \eqref{beka} holds with
    $$
    \ka=\min\Big\{\ka_3(j-1;\al,d),\,\,\ka_4(j;\al,d)\Big\}\,.
    $$
    It is straightforward to check that
    $$
    \ka_4(j;\al,d)\leq \ka_3(j-1;\al,d) \iff \al \leq d-j+\frac{d-2j+1}{d-j}\,,
    $$
    and
    $$
    0<\frac{d-2j+1}{d-j}<1\,.
    $$
    
    \item For $\al \in [d-j,d-j+1]$ with $j= \floor{\frac{d+1}{3}}+1$, by \eqref{a1} we can take 
    $$\ka=\ka_3\left(\floor{\frac{d+1}{3}};\al,d\right)\,.$$
    
    \item For $\al \in [d-j,d-j+1]$ with $\floor{\frac{d+1}{3}}+2 \leq j \leq \floor{\frac d2}$, by applying \eqref{a1} when $j=\floor{\frac{d+1}{3}}+2$ and applying \eqref{b1} otherwise we can take $\ka=\ka_3(j-2;\al,d)$, by \eqref{b2} we can take $\ka=\max\big\{\ka_3(j-1;\al,d),\,\,\ka_5(j-1;\al,d)\big\}$, and by \eqref{b3} we can take $\ka=\ka_5(j;\al,d)$. Therefore, \eqref{beka} holds if we choose $\ka$ to be
   $$
        \min\Big\{\ka_3(j-2;\al,d),\,\,\max\big\{\ka_3(j-1;\al,d),\,\,\ka_5(j-1;\al,d)\big\},\,\,\ka_5(j;\al,d)\Big\},
   $$
   and \eqref{b2} tells us that this number is 
   $$
   \min\Big\{\ka_3(j-1;\al,d),\,\,\ka_5(j;\al,d)\Big\} \quad \text{for} \quad \al \leq d-j+\frac{2(d-2j+1)}{d-j-2}
   $$
   and
   $$
   \min\Big\{\ka_3(j-2;\al,d),\,\,\ka_5(j-1;\al,d)\Big\} \quad \text{for} \quad \al \geq d-j+\frac{2(d-2j+1)}{d-j-2}.
   $$
   
   \item For $d$ odd, $d\geq 7$ and $\al\in(\frac{d-1}{2},\frac {d+1}{2}]$, by applying \eqref{a1} when $d=7,9,11$ and applying \eqref{b1} when $d\geq 13$, we can take $$\ka=\ka_3\left(\frac{d-3}{2};\al,d\right)\,.$$
   Note that when $d=3,5$, the case $\al\in(\frac{d-1}{2},\frac {d+1}{2}]$ is the same as the case $\al\in[d-j,d-j+1]$ with $j= \floor{\frac{d+1}{3}}+1$, and we have
 $$
\ka=\ka_3(1;\al,3) \text{ for } d=3, \quad \text{and} \quad  \ka=\ka_3(2;\al,5) \text{ for } d=5\,.
$$
   
   \item For $d$ even and $\al\in(\frac {d-1}{2},\frac {d}{2}]$, by applying \eqref{a1} when $d=4,6,8$ and applying \eqref{b1} when $d\geq 10$, we can take
   $$
   \ka=\ka_3\left(\frac{d}{2}-1;\al,d\right).
   $$
\end{itemize}

Note that the above discussion covers all the cases $d\geq 3$ and $\al\in\left(\frac{d-1}{2},d\right)$ for Theorem \ref{thm:prb}.
It remains to verify Lemmas \ref{C:sizeEf} and \ref{C:si}, and we will do so in the following two subsections.

\subsection{Proof of Lemma \ref{C:sizeEf}}
Since $f=\hichi_\Om$, we have
$$
Ef(Rx)=\frac{1}{(2\pi)^{d/2}} \int_{\Om} e^{iR(x'\cdot\xi'+x''\cdot \xi'' +x_d|\xi'|^2+x_d|\xi''|^2)} \,d\xi. 
$$
So it suffices to prove that 
\begin{equation} \label{phase}
 R(x'\cdot\xi'+x''\cdot \xi'' +x_d|\xi'|^2+x_d|\xi''|^2) \in 2\pi \ZZ +(-\frac{1}{100},\frac{1}{100})\,,
\end{equation}
provided that $\xi\in \Om$ and $x\in\La$. Indeed, by definitions of $\Om$ and $\La$, we write
$$
|\xi'|<cR^{-\frac 12}, \quad |x'|<cR^{-\frac 1 2}
$$
$$
\xi''=2\pi R^{-\ka}m+v,\quad \text{where} \quad m\in\ZZ^{d-m-1}, |m|< \frac{1}{2\pi}R^{\ka}, |v|<cR^{-1},
$$
$$
x''=R^{\ka-1}\ell+u,\quad \text{where} \quad \ell\in\ZZ^{d-m-1}, |\ell|< R^{1-\ka}, |u|<cR^{-1},
$$
and 
$$
x_d = \frac{1}{2\pi}R^{2\ka-1}k+\e, \quad \text{where} \quad k\in \ZZ, |k|< 2\pi R^{1-2\ka}, |\e|<cR^{-1}\,.
$$
Let us look at the four components in \eqref{phase} separately:
\begin{itemize}
    
    \item $|Rx'\cdot\xi'|<RcR^{-1/2}cR^{-1/2}=c^2$\,,
    \vspace{0.1 in}
    \item Since $|x_d|<1$,
    $$\left|Rx_d|\xi'|^2\right|<R c^2R^{-1}=c^2\,,$$

    \item For $Rx''\cdot\xi''$, we have
    \begin{align*}
        Rx''\cdot\xi''&=R(R^{\ka-1}\ell+u)\cdot(2\pi R^{-\ka}m+v)\\
        &= 2\pi\ell\cdot m+R^\ka \ell \cdot v+2\pi R^{1-\ka}u\cdot m +Ru\cdot v\,,
    \end{align*}
    where $2\pi\ell\cdot m\in 2\pi \ZZ$ and the other three terms are bounded by
    $$
    R^\ka R^{1-\ka}cR^{-1}+R^{1-\ka}cR^{-1}R^\ka+RcR^{-1}cR^{-1}=c+c+c^2R^{-1}.
    $$

    \item For $Rx_d|\xi''|^2$, we have
    \begin{align*}
        Rx_d|\xi''|^2&=R(\frac{1}{2\pi}R^{2\ka-1}k+\e)(2\pi R^{-\ka}m+v)\cdot(2\pi R^{-\ka}m+v)\\
        &= 2\pi k|m|^2+2 R^\ka k(m\cdot v)+\frac{1}{2\pi}R^{2\ka}k|v|^2\\
        &\quad + 4\pi^2\e R^{1-2\ka}|m|^2+4\pi \e R^{1-\ka}(m\cdot v)+\e R |v|^2\,,    \end{align*}
        where $ 2\pi k|m|^2\in 2\pi \ZZ$ and the other five terms are bounded by 
        \begin{align*}
            &R^\ka R^{1-2\ka}R^\ka cR^{-1}+R^{2\ka} 2\pi R^{1-2\ka} c^2 R^{-2}+cR^{-1}R^{1-2\ka}R^{2\ka}\\
            &+cR^{-1}R^{1-\ka}R^\ka cR^{-1}+cR^{-1}Rc^2R^{-2}\\
            =&c+2\pi c^2 R^{-1}+c+c^2R^{-1}+c^3R^{-2}\,.
        \end{align*}
\end{itemize}
Therefore, \eqref{phase} follows by taking $c$ sufficiently small, say $c=1/1000$.

\subsection{Proof of Lemma \ref{C:si}}
Recall that $d\mu=\hichi_\La\,dx$ and $\La$ is defined by
\begin{equation}\label{La'}
    \La:=\left[B^m_{cR^{-1/2}}\times \left(R^{\ka-1}\ZZ^{d-m-1}+B^{d-m-1}_{cR^{-1}}\right)\times\left(\frac{1}{2\pi}R^{2\ka-1}\ZZ+I_{cR^{-1}}\right)\right] \cap B^d(0,1)\,.
\end{equation}
We aim to prove that
$$
c_\al(\mu) \sim R^{\al-d}\,,
$$
by taking $\ka$ as stated in Lemma \ref{C:si}.
   
Recall that
$$
c_\al(\mu):=\sup_{x\in\ZR^d, r>0} \frac{\mu(B(x,r))}{r^\al}=\sup_{r>0}c_\al(\mu,r)\,,
$$
where
$$
c_\al(\mu,r):=\sup_{x\in\ZR^d} \frac{\mu(B(x,r))}{r^\al}\,.
$$
We will calculate $c_\al(\mu,r)$ directly from \eqref{La'}.
The important scales for $r$ are 
$R^{-1}, R^{\ka-1}, R^{2\ka-1}$ and $R^{-1/2}$. To compare the scales $R^{2\ka-1}$ and $R^{-1/2}$, we consider the two cases $\ka\leq 1/4$ and $\ka>1/4$ separately.

\vspace{.25in}
\textbf{\emph{Case I}: $\ka \leq\frac 14\,.$ } In this case, the important scales for $r$ are ordered as follows:
$$
R^{-1}<R^{\ka-1}<R^{2\ka-1}\leq R^{-1/2}<1\,.
$$
Now we calculate $C_\al(\mu,r)$ for different values of $r$.
\begin{itemize}
    \item For $0<r\leq R^{-1}$,
    $$
    c_\al(\mu,r)\sim \frac{r^d}{r^\al}=r^{d-\al}\,.
    $$
    
    \noindent Since $d-\al>0$, we have
    \begin{equation}\label{1}
        \sup_{0<r\leq R^{-1}} c_\al(\mu,r)\sim  c_\al(\mu,R^{-1})\sim R^{\al-d}\,.
    \end{equation}
    
    \item For $R^{-1}\leq r\leq R^{\ka-1}$,
    $$
    c_\al(\mu,r)\sim \frac{r^m\cdot R^{-(d-m)}}{r^\al}=r^{m-\al}R^{-(d-m)}\,.
    $$
    
    \noindent If $\al\leq m$, we have
     \begin{equation}\label{2}
        \sup_{R^{-1}\leq r\leq R^{\ka-1}} c_\al(\mu,r)\sim  c_\al(\mu,R^{\ka-1})\,,
    \end{equation}
    
    \noindent and if $\al\geq m$, we have
     \begin{equation}\label{3}
        \sup_{R^{-1}\leq r\leq R^{\ka-1}} c_\al(\mu,r)\sim  c_\al(\mu,R^{-1})\,.
    \end{equation}
    
    \item For $R^{\ka-1}\leq r\leq R^{2\ka-1}$,
    $$
    c_\al(\mu,r)\sim \frac{r^m\cdot \left(\frac{r}{R^{\ka-1}}R^{-1}\right)^{d-m-1}\cdot R^{-1}}{r^\al}=r^{d-1-\al}R^{-\ka (d-m-1)-1}\,.
    $$
   
    \noindent If $\al\leq d-1$, we have
    \begin{equation}\label{4}
        \sup_{R^{\ka-1}\leq r\leq R^{2\ka-1}} c_\al(\mu,r)\sim  c_\al(\mu,R^{2\ka-1})\,,
    \end{equation}
    
    \noindent and if $\al\geq d-1$, we have
    \begin{equation}\label{5}
        \sup_{R^{\ka-1}\leq r\leq R^{2\ka-1}} c_\al(\mu,r)\sim  c_\al(\mu,R^{\ka-1})\,.
    \end{equation}
    
    \item For $R^{2\ka-1}\leq r\leq R^{-\frac 1 2}$,
    $$
    c_\al(\mu,r)\sim \frac{r^m\cdot \left(\frac{r}{R^{\ka-1}}R^{-1}\right)^{d-m-1}\left(\frac{r}{R^{2\ka-1}}R^{-1}\right)}{r^\al}=r^{d-\al}R^{-\ka( d-m+1)}\,.
    $$
    
    \noindent Since $d-\al>0$, we have
    \begin{equation}\label{6}
        \sup_{R^{2\ka-1}\leq r\leq R^{-1/2}} c_\al(\mu,r)\sim  c_\al(\mu,R^{-1/2})\sim R^{-\frac{d-\al}{2}-\ka(d-m+1)}\,.
    \end{equation}
    
     \item For $R^{-\frac 1 2}\leq r\leq 1$,
    $$
    c_\al(\mu,r)\sim \frac{R^{-m/2}\cdot \left(\frac{r}{R^{\ka-1}}R^{-1}\right)^{d-m-1}\left(\frac{r}{R^{2\ka-1}}R^{-1}\right)}{r^\al}=r^{d-m-\al}R^{-\ka (d-m+1)-\frac m2}\,.
    $$
     
    \noindent If $\al\leq d-m$, we have
    \begin{equation}\label{7}
        \sup_{R^{-1/2}\leq r\leq 1} c_\al(\mu,r)\sim  c_\al(\mu,1)\sim R^{-\ka(d-m+1) -\frac m 2}\,,
    \end{equation}
    
    \noindent and if $\al\geq d-m$, we have
    \begin{equation}\label{8}
        \sup_{R^{-1/2}\leq r\leq 1} c_\al(\mu,r)\sim  c_\al(\mu,R^{-\frac 12})\,.
    \end{equation}
\end{itemize}

\noindent It is also obvious that 
$$
\sup_{r\geq 1} c_\al(\mu,r)\sim  c_\al(\mu,1)\,.
$$
Therefore, for $m\leq \al\leq d-m$, by combining \eqref{1}, \eqref{3}, \eqref{4}, \eqref{6} and \eqref{7}, we can tell that
\begin{align*}
    c_\al(\mu)&\sim \max\left\{c_\al(\mu, R^{-1}),c_\al(\mu,1)\right\}\\
    &\sim \max\left\{R^{\al-d}, R^{-\ka (d-m+1)-\frac m2}\right\} =R^{\al-d}\,,
\end{align*}
provided that
$$
\ka=\ka_3(m;\al,d)=\frac{d-m/2-\al}{d-m+1}\,.
$$

\noindent For $d-m\leq\al\leq d-1$, by combining \eqref{1}, \eqref{3}, \eqref{4}, \eqref{6} and \eqref{8}, we can tell that
\begin{align*}
    c_\al(\mu)&\sim \max\left\{c_\al(\mu, R^{-1}),c_\al(\mu,R^{-\frac 12})\right\}\\
    &\sim \max\left\{R^{\al-d}, R^{-\frac{d-\al}{2}-\ka(d-m+1)}\right\} =R^{\al-d}\,,
\end{align*}
provided that
$$
\ka=\ka_4(m;\al,d)=\frac{d-\al}{2(d-m+1)}\,.
$$

\noindent For $d-1\leq\al<d$, by combining \eqref{1}, \eqref{3}, \eqref{5}, \eqref{6} and \eqref{8}, we can tell that
$$
    c_\al(\mu) \sim \max\left\{c_\al(\mu, R^{-1}),c_\al(\mu,R^{-\frac 12})\right\}\sim R^{\al-d}\,,
$$
provided that
$$
\ka=\ka_4(m;\al,d)\,.
$$

\noindent Note that the calculation of $c_\al(\mu)$ above is in the case $\ka\leq 1/4$. While
$$
\ka_3(m;\al,d) \leq \frac 14 \iff \al \geq \frac{3d-m-1}{4}\,,
$$
and
$$
\ka_4(m;\al,d) \leq \frac 14 \iff \al \geq \frac{d+m-1}{2}\,.
$$
Also note that
$$
m < \frac{3d-m-1}{4} \quad \text{ for } \quad d\geq 3\,,
$$
$$
\frac{3d-m-1}{4} \leq d-m \iff m\leq \frac{d+1}{3}\,,
$$
and
$$
d-m \geq \frac{d+m-1}{2} \iff m\leq \frac{d+1}{3}\,.
$$
Therefore, in \emph{Case I} we obtain $c_\al(\mu)\sim R^{\al-d}$ by taking $\ka$ as follows:
\begin{itemize}
    \item If $1\leq m\leq \frac{d+1}{3}$, then
    \begin{equation} \label{A1}
    \ka=\ka_3(m;\al,d)\quad \text{for } \quad \frac{3d-m-1}{4}\leq \al\leq d-m, 
\end{equation}
and
\begin{equation} \label{A2}
    \ka=\ka_4(m;\al,d)\quad \text{for } \quad d-m\leq \al< d. 
\end{equation}

\item If $\frac{d+1}{3} < m\leq \frac d 2$, then
\begin{equation}\label{A3}
    \ka=\ka_4(m;\al,d)\quad \text{for } \quad  \frac{d+m-1}{2} \leq \al <d.
\end{equation}

\end{itemize}

\vspace{.25in}

\textbf{\emph{Case II}: $\ka >\frac 14\,.$ } Note that, we have proved Lemma \ref{C:si} for $\al\geq d-1$ in \emph{Case I}. Therefore, here we can assume that $\al<d-1$. In this case, the important scales for $r$ are ordered as follows:
$$
R^{-1}<R^{\ka-1}< R^{-1/2}< R^{2\ka-1}<1\,.
$$
Now we calculate $C_\al(\mu,r)$ for different values of $r$.
\begin{itemize}
    \item For $0<r\leq R^{\ka-1}$, same as in \emph{Case I}, if $\al\leq m$ we have
    \begin{equation}\label{1'}
        \sup_{0<r\leq R^{\ka-1}} c_\al(\mu,r)\sim  c_\al(\mu,R^{\ka-1})\,,
    \end{equation}
    
    \noindent and if $\al\geq m$ we have
    \begin{equation}\label{2'}
        \sup_{0<r\leq R^{\ka-1}} c_\al(\mu,r)\sim  c_\al(\mu,R^{-1})\sim R^{\al-d}\,.
    \end{equation}
    
    \item For $R^{\ka-1}\leq r\leq R^{-\frac 12}$,
    $$
    c_\al(\mu,r)\sim \frac{r^m\cdot \left(\frac{r}{R^{\ka-1}}R^{-1}\right)^{d-m-1}\cdot R^{-1}}{r^\al}=r^{d-1-\al}R^{-\ka (d-m-1)-1}\,.
    $$
    
    \noindent Since $\al<d-1$, we have
    \begin{equation}\label{3'}
        \sup_{R^{\ka-1}\leq r\leq R^{-1/2}} c_\al(\mu,r)\sim  c_\al(\mu,R^{-\frac 12})\sim R^{-\frac{d-\al}{2}-\ka (d-m-1)-\frac 12}\,.
    \end{equation}

    \item For $R^{-\frac 1 2}\leq r\leq R^{2\ka-1}$,
    $$
    c_\al(\mu,r)\sim \frac{R^{-\frac m 2}\cdot \left(\frac{r}{R^{\ka-1}}R^{-1}\right)^{d-m-1} \cdot R^{-1}}{r^\al}=r^{d-m-1-\al}R^{-\ka (d-m-1)-\frac m2-1}\,.
    $$
   
   \noindent If $\al\leq d-m-1$, we have
    \begin{equation}\label{4'}
        \sup_{R^{-1/2}\leq r\leq R^{2\ka-1}} c_\al(\mu,r)\sim  c_\al(\mu,R^{2\ka-1})\,,
    \end{equation}
    
    \noindent and if $\al\geq d-m-1$, we have
    \begin{equation}\label{5'}
        \sup_{R^{- 1/2}\leq r\leq R^{2\ka-1}} c_\al(\mu,r)\sim  c_\al(\mu,R^{-\frac 12})\,.
    \end{equation}
    
     \item For $R^{2\ka-1}\leq r\leq 1$,
    $$
    c_\al(\mu,r)\sim \frac{R^{-\frac m 2}\cdot \left(\frac{r}{R^{\ka-1}}R^{-1}\right)^{d-m-1}\left(\frac{r}{R^{2\ka-1}}R^{-1}\right)}{r^\al}=r^{d-m-\al}R^{-\ka (d-m+1)-\frac m2}\,.
    $$
    
    \noindent If $\al\leq d-m$, we have
    \begin{equation}\label{6'}
        \sup_{R^{2\ka-1}\leq r\leq 1} c_\al(\mu,r)\sim  c_\al(\mu,1)\sim R^{-\ka (d-m+1) -\frac m 2}\,.
    \end{equation}
    
    \noindent and if $\al\geq d-m$, we have
    \begin{equation}\label{7'}
        \sup_{R^{2\ka-1}\leq r\leq 1} c_\al(\mu,r)\sim  c_\al(\mu,R^{2\ka-1})\,.
    \end{equation}
\end{itemize}

\noindent It is also obvious that
$$
        \sup_{r\geq 1} c_\al(\mu,r)\sim  c_\al(\mu,1)\,.
$$
Therefore, for $m\leq \al\leq d-m-1$, by combining \eqref{2'},  \eqref{3'}, \eqref{4'} and \eqref{6'}, we can tell that
\begin{align*}
    c_\al(\mu)&\sim \max\left\{c_\al(\mu, R^{-1}),c_\al(\mu,1)\right\}\\
    &\sim \max\left\{R^{\al-d}, R^{-\ka (d-m+1)-\frac m2}\right\} =R^{\al-d}\,,
\end{align*}
provided that
$$
\ka=\ka_3(m;\al,d)=\frac{d-m/2-\al}{d-m+1}\,.
$$

For $m\leq d-m-1 \leq \al\leq d-m$ (and so $m<d/2$), by combining \eqref{2'}, \eqref{3'}, \eqref{5'} and \eqref{6'}, we can tell that
\begin{align*}
    c_\al(\mu)&\sim \max\left\{c_\al(\mu, R^{-1}),c_\al(\mu,R^{-\frac 12}), c_\al(\mu,1)\right\}\\
    &\sim \max\left\{R^{\al-d}, R^{-\frac{d-\al}{2}-\ka (d-m-1)-\frac 12}, R^{-\ka (d-m+1)-\frac m2}\right\} =R^{\al-d}\,,
\end{align*}
provided that
$$
\ka\geq \ka_3(m;\al,d) \quad \text{and} \quad \ka\geq \ka_5(m;\al,d)=\frac{d-\al-1}{2(d-m-1)}\,.
$$
Therefore, we can take 
$$
\ka=\max\left\{ \ka_3(m;\al,d), \ka_5(m;\al,d) \right\}\,.
$$
While, by a direct calculation, if $m\leq \frac{d+1}{3}$, then 
$$
\ka=\ka_3(m;\al,d), \quad \text{for} \quad d-m-1\leq \al\leq d-m\,;
$$
and if $m> \frac{d+1}{3}$, then
$$
    \ka=\begin{cases}
    \ka_3(m;\al,d), \quad &\text{for} \quad d-m-1\leq\al\leq d-m-1+\frac{2(d-2m-1)}{d-m-3}\,, \vspace{0.1in} \\
    \ka_5(m;\al,d), \quad &\text{for}\quad  d-m-1+\frac{2(d-2m-1)}{d-m-3}\leq \al\leq d-m\,.
    \end{cases}
$$

Next, for $m\leq d-m\leq \al< d-1$ (and so $m\leq d/2$), by combining \eqref{2'}, \eqref{3'}, \eqref{5'} and \eqref{7'}, we can tell that
$$
    c_\al(\mu)\sim \max\left\{c_\al(\mu, R^{-1}),c_\al(\mu,R^{-\frac 12})\right\} \sim R^{\al-d}\,,
$$
provided that
$$
\ka=\ka_5(m;\al,d)\,.
$$

Note that the calculation of $c_\al(\mu)$ above is in the case $\ka\geq 1/4$. While
$$
\ka_3(m;\al,d) \geq \frac 14 \iff \al \leq \frac{3d-m-1}{4}\,,
$$
and
$$
\ka_5(m;\al,d) \geq \frac 14 \iff \al \leq \frac{d+m-1}{2}\,.
$$
Also, note that 
$$
\frac{3d-m-1}{4} \geq d-m-1  \iff m\geq \frac{d-3}{3}\,, 
$$
$$
\frac{3d-m-1}{4} \geq d-m  \iff m\geq \frac{d+1}{3}\,, 
$$
and 
$$
\frac{d+m-1}{2} \geq d-m  \iff m\geq \frac{d+1}{3}\,, 
$$
Therefore, in \emph{Case II} we obtain $c_\al(\mu)\sim R^{\al-d}$ by taking $\ka$ as follows:
\begin{itemize}
    \item If $1\leq m\leq \frac{d+1}{3}$, then
    \begin{equation} \label{B1}
    \ka=\ka_3(m;\al,d)\quad \text{for } \quad m\leq \al\leq \frac{3d-m-1}{4}.
\end{equation}

\item If $\frac{d+1}{3} < m < \frac d 2$, then
\begin{equation}\label{B2}
    \ka=\ka_3(m;\al,d)\quad \text{for } \quad m\leq \al\leq d-m-1, 
\end{equation}
and
\begin{align} \label{B3}
    \ka&=\max\Big\{\ka_3(m;\al,d),\,\,\ka_5(m;\al,d)\Big\} \quad \text{for } \quad d-m-1\leq \al\leq d-m \\
    &=
    \begin{cases}
    \ka_3(m;\al,d), & d-m-1 \leq \al \leq d-m-1+\frac{2(d-2m-1)}{d-m-3}\,,\vspace{0.1 in} \notag\\
    \ka_5(m;\al,d), & d-m-1+\frac{2(d-2m-1)}{d-m-3}\leq \al \leq d-m \,,
    \end{cases}
\end{align}
and
\begin{equation}\label{B4}
    \ka=\ka_5(m;\al,d)\quad \text{for } \quad  d-m\leq \al\leq \frac{d+m-1}{2}.
\end{equation}
And \eqref{B4} also holds and is nontrivial when $m=\frac d2$\,.

\end{itemize}

The proof of Lemma \ref{C:si} is done by combining the conclusions from both \emph{Case I} and \emph{Case II}.

\end{document}